# Decomposition of Almost Poisson Structure of Non-Self-Adjoint Dynamical Systems

Yongxin Guo, Chang Liu, Shixing Liu and Peng Chang

College of Physics, Liaoning University

Shenyang 110036, China

Email: yxguolnu.edu.cn

**Abstract.** Non-self-adjoint dynamical systems, *e.g.*, nonholonomic systems, can admit an almost Poisson structure, which is formulated by a kind of Poisson bracket satisfying the usual properties except for the Jacobi identity. A general theory of the almost Poisson structure is investigated based on a decomposition of the bracket into a sum of a Poisson one and an almost Poisson one. The corresponding relation between Poisson structure and symplectic structure is proved, making use of Jacobiizer and symplecticizer. Based on analysis of pseudo-symplectic structure of constraint submanifold of Chaplygin's nonholonomic systems, an almost Poisson bracket for the systems is constructed and decomposed into a sum of a canonical Poisson one and an almost Poisson one. Similarly, an almost Poisson structure, which can be decomposed into a sum of canonical one and an almost "Lie-Poisson" one, is also constructed on an affine space with torsion whose autoparallels are utilized to described the free motion of some non-self-adjoint systems. The decomposition of the almost Poisson bracket directly leads to a decomposition of a dynamical vector field into a sum of usual Hamiltonian vector field and an almost Hamiltonian one, which is useful to simplifying the integration of vector fields.

## 1. Introduction

In the framework of inverse problem of dynamics, dynamical systems can be classified into self-adjoint and non-self-adjoint ones. From the viewpoint of calculus of variation, a system of ordinary variational forms is termed self-adjoint when it coincides with its adjoint system for all admissible variations. A set of ordinary differential equations is called self-adjoint if the corresponding variational forms are self-adjoint. Otherwise it is called non-self-adjoint[1]. In themodern setting of differential geometry, the self-adjointness of the dynamical systems can also be equivalently defined by the conditions satisfied by symmetries of equations of motion. A dynamical system is called self-adjoint if its dynamical symmetries coincide with its adjoint symmetries[2,3]. Evidently conservative systems are self-adjoint. The converse is not true. For the Newtonian systems in the fundamental form or kinematic form, the conditions of self-adjointness of the systems are Helmholtz's conditions, which lead to a direct Lagrangian or Hamiltonian representation of the systems. Geometrically, self-adjointness of Lagrangian or Hamiltonian systems can be proved to be accordant with symplecticity of phase space. So such self-adjoint systems can admit Poisson structure and easy to integrate.

Poisson structure for self-adjoint dynamical systems is formulated by Poisson brackets on the set of functions on manifold with the property of anticommutativity, bilinearity, Leibniz's rule and Jacobi's identity[4]. As well known, Lagrangian or Hamiltonian representation, either direct or indirect by self-adjoint genotopic transformation[5], in the local coordinates and time variables actually used in experiments is not universal. Universality of Lagrangian or Hamiltonian representation is only indirect in the sense of using Darboux's transformation of symplectic geometry. Therefore, many dynamical systems can not universally admit a direct Poisson structure, especially for the essentially non-self-adjoint dynamical systems. Even the direct Poisson structure does not exist for non-autonomous Birkhoffian systems and generalized Birkhoffian systems (nonlocal non-self-adjoint systems)[5].

There exist many non-self-adjoint physical systems such as the network modelling of energy conserving physical systems with external ports[6], nonholonomic constrained systems[7-11], some physical structure closely related with torsion of general affine metric space or spacetime *e. g.*, particles moving in Riemann-Cartan spacetime[12-18], a crystal with dislocation[19], motion of rigid body in body-fixed coordinate system[20], *etc*. Their configuration or phase space can admit an almost (quasi- or pseudo-) Poisson structure[21-34] in the sense that a kind of bracket existing on the set of functions on the manifold shares the usual properties of a Poisson bracket except for the Jacobi's identity. The equation of motion of the non-self-adjoint systems with the almost Poisson structure is much more difficult to resolve than that of self-adjoint systems with Poisson structure. However, in many cases the almost Poisson structure can be simplified by means of a decomposition of the bracket into a sum of canonical Poisson one and an almost Lie-Poisson one. In this article, we give a technique of decomposition of almost Poisson bracket and the corresponding dynamical vector for some non-self-adjoint dynamical systems. In section 2 a general theory of almost Poisson structure is formulated based on the decomposition technique. The close relation between the Poisson structure and symplectic structure for even dimensional manifold is proved with the help of Jacobiizer and symplecticizer. In section 3 pseudo-symplectic structure on the constraint submanifold is constructed for Chaplygin's nonholonomic systems, which leads to an almost Poisson structure by Legendre transformation. This almost Poisson structure is proved to be decomposed into a sum of canonical Poisson one and an almost Poisson one depending on the nonholonomicity of constraints. In section 4 an almost Poisson structure is similarly constructed on affine space with torsion, whose autoparallels deviate from its geodesics and is utilized to formulate the motion of many non-self-adjoint dynamical systems. Based on an analysis of inverse problem of calculus of variations for the autoparallels, an almost Poisson structure is constructed, depending on the torsion of the space. Such an almost Poisson structure can also be decomposed into a sum of canonical Poisson one and an almost Lie-Poisson one relating with torsion tensor of the affine space. The Einstein's summation convention is used throughout the article.

## 2. Decomposition of the almost Poisson structure on a manifold

Let $M$ be a dimensional manifold and let $\mathcal{F}(M)$ denote the set of smooth real-valued functions on $M$. Define a bracket operation denoted by $[\ ,\ ]: \mathcal{F}(M) \times \mathcal{F}(M) \to \mathcal{F}(M)$ satisfying the following relation:

$$[f, g] = -[g, f] \tag{2.1a}$$

$$[c_1 f + c_2 g, h] = c_1[f, h] + c_2[g, h], \quad c_1, c_2 \in R \tag{2.1b}$$

$$[fg,h] = f[g,h] + [f,h]g \tag{2.1c}$$

This bracket is called almost Poisson bracket. Accordingly the pair $\{M,[\ ,\ ]\}$ is called almost Poisson manifold.

If the bracket further satisfies the Jacobi's identity:

$$[[f,g],h] + [[g,h],f] + [[h,f],g] = 0 \tag{2.2}$$

the $[\ ,\ ]$ becomes Poisson bracket and the manifold $M$ Poisson one. The conditions (2.1)-(2.2) make $(\mathcal{F}(M),[\ ,\ ])$ into a Lie algebra.

Suppose the almost Poisson bracket $[\ ,\ ]$ can be decomposed into a sum of the Poisson bracket $\{\ ,\ \}$ satisfying the conditions (2.1)-(2.2) and an bracket $[\![\ ,\ ]\!]$, i.e.,

$$[f,g] = \{f,g\} + [\![f,g]\!] \tag{2.3}$$

It can be proved that the new bracket $[\![\ ,\ ]\!]$ satisfies:

$$[\![f,g]\!] = -[\![g,f]\!] \tag{2.4a}$$

$$[\![c_1 f + c_2 g, h]\!] = c_1 [\![f,h]\!] + c_2 [\![g,h]\!], \quad c_1, c_2 \in R \tag{2.4b}$$

$$[\![fg,h]\!] = f[\![g,h]\!] + [\![f,h]\!]g \tag{2.4c}$$

i.e., it is an almost Poisson bracket. Hence, without proof we can obtain the following theorem:

**Theorem 2.1**: Any almost Poisson bracket can be decomposed into a sum of a Poisson bracket and another almost Poisson bracket on the same manifold. Conversely, a Poisson bracket plus an almost Poisson bracket can also be verified to lead to another almost Poisson bracket.

We should point that this relation of decomposition cannot generalize to the Poisson bracket. Substitute the Eq. (2.3) into left side of Eq. (2.2), considering the anticommutativity, bilinearity of the almost Poisson bracket, we get

$$\begin{aligned}
&[[f,g],h] + [[g,h],f] + [[h,f],g] \\
&= [\{f,g\} + [\![f,g]\!], h] + [\{g,h\} + [\![g,h]\!], f] + [\{h,f\} + [\![h,f]\!], g] \\
&= \{\{f,g\} + [\![f,g]\!], h\} + [\![\{f,g\} + [\![f,g]\!], h]\!] + \{\{g,h\} + [\![g,h]\!], f\} + [\![\{g,h\} + [\![g,h]\!], f]\!] \\
&\quad + \{\{h,f\} + [\![h,f]\!], g\} + [\![\{h,f\} + [\![h,f]\!], g]\!] \\
&= \{\{f,g\},h\} + \{\{g,h\},f\} + \{\{h,f\},g\} + [\![[\![f,g]\!],h]\!] + [\![[\![g,h]\!],f]\!] + [\![[\![h,f]\!],g]\!] \\
&\quad + \{[\![f,g]\!],h\} + \{[\![g,h]\!],f\} + \{[\![h,f]\!],g\} + [\![\{f,g\},h]\!] + [\![\{g,h\},f]\!] + [\![\{h,f\},g]\!]
\end{aligned} \tag{2.5}$$

Since the bracket $\{\ ,\ \}$ satisfies the Jacobi's identity:

$$\{\{f,g\},h\}+\{\{g,h\},f\}+\{\{h,f\},g\}=0 \tag{2.6}$$

then the Eq. (2.5) becomes

$$[[f,g],h]+[[g,h],f]+[[h,f],g]=[[[f,g],h]]+[[[g,h],f]]+[[[h,f],g]] \\ +\{[[f,g]],h\}+\{[[g,h]],f\}+\{[[h,f]],g\} \tag{2.7} \\ +[[\{f,g\},h]]+[[\{g,h\},f]]+[[\{h,f\},g]]$$

which shows that the brackets $[\,,\,]$ and $[\![\,,\,]\!]$ do not likely satisfies the Jacobi's identity simultaneously unless the following relation exists:

$$\{[[f,g]],h\}+\{[[g,h]],f\}+\{[[h,f]],g\}+[[\{f,g\},h]]+[[\{g,h\},f]]+[[\{h,f\},g]]=0 \tag{2.8}$$

Considering this and theorem 2.1, we obtain the following theorem:

**Theorem 2.2:** A Poisson bracket $[\,,\,]$ can be decomposed into a sum of two Poisson brackets $\{\,,\,\}$ and $[\![\,,\,]\!]$, i.e., $[f,g]=\{f,g\}+[\![f,g]\!]$ for $f,g\in\mathcal{F}(M)$ or a sum of a Poisson bracket $\{\,,\,\}$ and a Poisson bracket $[\![\,,\,]\!]$ is also a Poisson bracket if and only if such two Poisson brackets $\{\,,\,\}$ and $[\![\,,\,]\!]$ are coupled by the Eq. (2.8).

Therefore, unlike the almost Poisson bracket, a Poisson bracket can not be decomposed into a sum of two Poisson brackets in general and vice versa. Of course, a sum of two Poisson brackets is an almost Poisson bracket.

The anti-commutativity (2.1a), bilinearity (2.1b) and Leibniz's rule (2.1c) lead to the existence of an anti-symmetric tensor on $M$, denoted by $\mathbf{J}$ which assigns to each point $x\in M$ a linear map $\mathbf{J}(x):T_x^*M\to T_xM$ such that

$$[f,g]=\langle\mathbf{J}(x)\cdot\mathrm{d}f(x),\mathrm{d}g(x)\rangle \tag{2.9}$$

where $\langle\,,\,\rangle$ denotes the usual pair of vectors and forms on the manifold $M$. The tensor $\mathbf{J}$ can be called almost Poisson tensor on $M$. Let $\dim M=m$ and $\{x^i\}(i=1,2,\cdots,m)$ denotes local coordinates on $M$. Then the Eq. (2.9) can be represented in coordinates

$$[f,g]=\mathbf{J}^{ij}\frac{\partial f}{\partial x^i}\frac{\partial g}{\partial x^j} \tag{2.10}$$

where the summation convention is used. This equation obviously gives

$$[x^i,x^j]=\mathbf{J}^{ij} \tag{2.11}$$

According to Eq. (2.3) and Eq. (2.11), the almost Poisson tensor $\mathbf{J}^{ij}$ can be decomposed into

$$\mathbf{J}^{ij}=\boldsymbol{\omega}^{ij}+\mathbf{K}^{ij} \tag{2.12}$$

where $\boldsymbol{\omega}^{ij} = \{x^i, x^j\}$ is a Poisson tensor and $\mathbf{K}^{ij} = [\![x^i, x^j]\!]$ is an almost Poisson tensor. In terms of these tensors the Jacobi's identity, e. g., Eq. (2.2) becomes

$$\mathbf{J}^{lk}\frac{\partial \mathbf{J}^{ij}}{\partial x^l} + \mathbf{J}^{li}\frac{\partial \mathbf{J}^{jk}}{\partial x^l} + \mathbf{J}^{lj}\frac{\partial \mathbf{J}^{ki}}{\partial x^l} = 0 \tag{2.13}$$

The couple relation (2.8) of brackets $\{\,,\,\}$ and $[\![\,,\,]\!]$ becomes

$$\boldsymbol{\omega}^{lk}\frac{\partial \mathbf{K}^{ij}}{\partial x^l} + \boldsymbol{\omega}^{li}\frac{\partial \mathbf{K}^{jk}}{\partial x^l} + \boldsymbol{\omega}^{lj}\frac{\partial \mathbf{K}^{ki}}{\partial x^l} + \mathbf{K}^{lk}\frac{\partial \boldsymbol{\omega}^{ij}}{\partial x^l} + \mathbf{K}^{li}\frac{\partial \boldsymbol{\omega}^{jk}}{\partial x^l} + \mathbf{K}^{lj}\frac{\partial \boldsymbol{\omega}^{ki}}{\partial x^l} = 0 \tag{2.14}$$

Here we have utilized the relations

$$\{f,g\} = \boldsymbol{\omega}^{ij}\frac{\partial f}{\partial x^i}\frac{\partial g}{\partial x^j}, \qquad [\![f,g]\!] = \mathbf{K}^{ij}\frac{\partial f}{\partial x^i}\frac{\partial g}{\partial x^j} \tag{2.15}$$

The theorem 2.1 and theorem 2.2 can be applied to the case of tensor representation.

**Corollary 2.1:** Any almost Poisson tensor can be decomposed into a sum of a Poisson tensor and another almost Poisson tensor on the same manifold. Conversely, a Poisson tensor plus an almost Poisson tensor can also be verified to lead to another almost Poisson tensor.

**Corollary 2.2:** A Poisson tensor $\mathbf{J}^{ij}$ can be decomposed into a sum of two Poisson tensors $\boldsymbol{\omega}^{ij}$ and $\mathbf{K}^{ij}$, i.e., $\mathbf{J}^{ij} = \boldsymbol{\omega}^{ij} + \mathbf{K}^{ij}$ or a sum of a Poisson tensor $\boldsymbol{\omega}^{ij}$ and a Poisson tensor $\mathbf{K}^{ij}$ is also a Poisson tensor if and only if such two Poisson tensors $\boldsymbol{\omega}^{ij}$ and $\mathbf{K}^{ij}$ are coupled by the Eq. (2.14).

If $\dim M = m = 2n$ and the almost Poisson tensor $\mathbf{J}$ on $M$ is non-degenerate in the sense that it has maximal rank at each point $x \in M$, i.e., $\text{rank}\mathbf{J} = 2n$, there exist a fundamental 2-form $\boldsymbol{\Omega}$ on $M$ which is an inverse of the almost Poisson tensor $\mathbf{J}$. Such a fundamental 2-form is obviously non-degenerate. But it is a pseudo-symplectic form without use of Darboux's transformation. In this paper we discuss direct geometric structure such as almost Poisson structure, pseudo-symplectic structure in the sense of preserving the original coordinates on the almost Poisson manifold $M$. We can prove the following important theorem:

**Theorem 2.3:** The fundamental 2-form $\boldsymbol{\Omega} = \mathbf{J}^{-1}$ on $M$ is symplectic if and only if the tensor $\mathbf{J} = \boldsymbol{\Omega}^{-1}$ is Poisson.

Proof. Assume that $\boldsymbol{\Omega} = \mathbf{J}^{-1}$ is a symplectic form. $\boldsymbol{\Omega}$ can be represented in local coordinates by

$$\boldsymbol{\Omega} = \boldsymbol{\Omega}_{ij} dx^i \wedge dx^j \tag{2.16}$$

$d\boldsymbol{\Omega} = 0$ leads to

$$\frac{\partial \boldsymbol{\Omega}_{ij}}{\partial x^k} + \frac{\partial \boldsymbol{\Omega}_{jk}}{\partial x^i} + \frac{\partial \boldsymbol{\Omega}_{ki}}{\partial x^j} = 0 \tag{2.17}$$

Here we directly define an almost Poisson tensor $\mathbf{J} = \mathbf{\Omega}^{-1}$. Define a Jacobitizer and a symplecticizer respectively by

$$\mathbf{J}^{kij} \triangleq \mathbf{J}^{lk}\frac{\partial \mathbf{J}^{ij}}{\partial x^l} + \mathbf{J}^{li}\frac{\partial \mathbf{J}^{jk}}{\partial x^l} + \mathbf{J}^{lj}\frac{\partial \mathbf{J}^{ki}}{\partial x^l}, \quad \mathbf{\Omega}_{kij} \triangleq \frac{\partial \mathbf{\Omega}_{ij}}{\partial x^k} + \frac{\partial \mathbf{\Omega}_{jk}}{\partial x^i}\frac{\partial \mathbf{\Omega}_{ki}}{\partial x^j} \tag{2.18}$$

Making use of

$$\frac{\partial \mathbf{\Omega}_{ij}}{\partial x^k} = \mathbf{\Omega}_{jn}\frac{\partial \mathbf{J}^{nl}}{\partial x^k}\mathbf{\Omega}_{li} \tag{2.19}$$

we obtain

$$\mathbf{\Omega}_{kij}\mathbf{J}^{km}\mathbf{J}^{in}\mathbf{J}^{jp} = \mathbf{J}^{mnp} \tag{2.20}$$

Obviously the symplectic condition $\mathbf{\Omega}_{kij} = 0$ leads to Jacobi's identity $\mathbf{J}^{mnp} = 0$. Conversely, if a Poisson structure is given by a tensor $\mathbf{J}$ satisfying the Eq. (2.13), a fundamental 2-form $\mathbf{\Omega} = \mathbf{J}^{-1}$ can then be defined. With the help of

$$\frac{\partial \mathbf{J}^{ij}}{\partial x^k} = \mathbf{J}^{jn}\frac{\partial \mathbf{\Omega}_{nl}}{\partial x^k}\mathbf{J}^{li} \tag{2.21}$$

it can be verified that

$$\mathbf{J}^{mnp}\mathbf{\Omega}_{mk}\mathbf{\Omega}_{ni}\mathbf{\Omega}_{pj} = \mathbf{\Omega}_{kij} \tag{2.22}$$

which means that the symplectic condition $\mathbf{\Omega}_{kij} = 0$ can be derived directly from the self-adjoint condition $\mathbf{J}^{mnp} = 0$. The theorem obviously implies that almost Poisson structure is closely related with the pseudo-symplectic structure.

Define an almost Hamiltonian vector field $X_f$ of a function $f \in \mathcal{F}(M)$ on $M$ by

$$i_{X_f}\mathbf{\Omega} \triangleq \mathbf{\Omega}(X_f, \cdot) = df \tag{2.23}$$

Then the almost Poisson bracket can be defined by

$$\mathbf{\Omega}(X_f, X_g) = i_{X_g}df = X_g(f) = [f, g] \tag{2.24}$$

which satisfies the Eqs. (2.1a)-(2.1c). This bracket does not satisfy the Jacobi's identity because the fundamental 2-form is not closed. So the vector field $X_f$ defined by the Eq. (2.23) is an almost Hamiltonian vector field. From the Eq. (2.24) the almost Hamiltonian vector field can also be represented by

$$X_g^i = \dot{x}^i = \mathbf{J}^{ij}\frac{\partial g}{\partial x^j} \tag{2.25}$$

The almost Hamiltonian vector field $X_g$ can also be defined by $X_g(f) = [f, g]$ if the almost

Poisson structure is given beforehand. Similarly a Hamiltonian vector field $\mathcal{X}_g$ and an almost Hamiltonian one can be defined by

$$\mathcal{X}_g(f) = \{f, g\}, \quad \mathbb{X}_g(f) = [\![f, g]\!] \tag{2.26}$$

with components

$$\mathcal{X}_g^i = \dot{x}_H^j = \boldsymbol{\omega}^{ij} \frac{\partial g}{\partial x^j}, \quad \mathbb{X}_g^i = \dot{x}_{nH}^j = \mathbf{K}^{ij} \frac{\partial g}{\partial x^j} \tag{2.27}$$

Thus by means of the decomposition relation (2.3) and definition (2.24) we get a decomposition relation of the almost Hamiltonian vector field:

$$X_g = \mathcal{X}_g + \mathbb{X}_g \tag{2.28}$$

whose components are

$$X_g^i = \mathcal{X}_g^i + \mathbb{X}_g^i, \quad \text{or} \quad \dot{x}_g^j = \boldsymbol{\omega}^{ij} \frac{\partial g}{\partial x^j} + \mathbf{K}^{ij} \frac{\partial g}{\partial x^j} \tag{2.29}$$

It should be pointed that the vector fields $X_g$, $\mathcal{X}_g$ and $\mathbb{X}_g$ are associated with $[\,,\,]$, $\{,\}$ and $[\![\,,\,]\!]$ respectively. A vector field $X_g$ is Hamiltonian if the 2-form $\boldsymbol{\Omega}$ is symplectic or the bracket $[\,,\,]$ is Poisson. According to the definition (2.24) and (2.26), theorem 2.1 means that an almost Hamiltonian vector can be decomposed into a sum of a Hamiltonian vector field and another almost Hamiltonian one. The application of theorem 2.2 to the almost Hamiltonian vector field leads to

**Corollary 2.3:** For functions $f, g, h \in F(M)$, a Hamiltonian vector field $X_g$ can be decomposed into a sum of two Hamiltonian vector fields $\mathcal{X}_g$ and $\mathbb{X}_g$, i.e., $X_g = \mathcal{X}_g + \mathbb{X}_g$ or a sum of one Hamiltonian vector field $\mathcal{X}_g$ and another Hamiltonian vector field $\mathbb{X}_g$ is also a Hamiltonian vector field if and only if such two Hamiltonian vector fields $\mathcal{X}_g$ and $\mathbb{X}_g$ satisfy

$$\mathcal{X}_h(\mathbb{X}_g(f)) + \mathcal{X}_f(\mathbb{X}_h(g)) + \mathcal{X}_g(\mathbb{X}_f(h)) + \mathbb{X}_h(\mathcal{X}_g(f)) + \mathbb{X}_f(\mathcal{X}_h(g)) + \mathbb{X}_g(\mathcal{X}_f(h)) = 0 \tag{2.30}$$

Finally, we discuss a special case of above decomposition of almost Poisson bracket that the Poisson bracket $\{,\}$ is canonical, i.e., $\boldsymbol{\omega}^{ij}$ is a canonical Poisson tensor or the symplectic form $\boldsymbol{\omega}_{ij}$ is simple. In this case the coordinates $\{x^i\}$ ($i = 1, 2, \cdots, 2n$) on $M$ can be classified into $\{q^s, p_s\}$ ($s = 1, 2, \cdots, n$), which are conjugate with each other with respect to $\boldsymbol{\omega}^{ij}$. The fundamental 2-form becomes $\boldsymbol{\Omega} = \mathrm{d}p_s \wedge \mathrm{d}q^s + \mathbf{K}_{ij} \mathrm{d}x^i \wedge \mathrm{d}x^j$. Furthermore, sometimes we can assume that $\mathbf{K}^{ij}$ or $[\![,]\!]$ is Poisson, e.g., Lie-Poisson one. Then the almost Hamiltonian vector field is decomposed into a canonical Hamiltonian vector field and a non-canonical Hamiltonian one, which may be easy to integrate respectively.

## 3. Almost Poisson structure of Chaplygin's nonholonomic systems and its decomposition

We consider a mechanical system constrained by linear nonholonomic constraints, called Chaplygin's nonholonomic system. Denote configuration manifold by $Q$ with local coordinates $\{q^s\}$ ($i=1,2,\cdots,n$). Let $TQ$ and $T^*Q$ be tangent bundle and cotangent bundle to $Q$ respectively. The Lagrangian of the system is denoted by $L(q,\dot{q}) \in C^2(TQ)$. The constraints is

$$\dot{q}^\alpha = B_\mu^\alpha(q^\nu)\dot{q}^\mu \quad (\mu=1,2,\cdots,n\text{-}g; \alpha=1,2,\cdots,g) \tag{3.1}$$

This construction of constraints distinguishes two sets of coordinates $\{q^\mu\}$ and $\{q^\alpha\}$, which makes the configuration manifold $Q$ be of fibred structure over its submanifold $Q_0$ with local coordinates $\{q^\mu\}$. Let the map $\pi: Q \to Q_0$ be a projection with respect to which the vertical space $VQ$ is obviously the kernel of differential mapping $d\pi$. The Ehresmann connection $h$ defined by the splitting of the exact sequence of $0 \to VQ \to TQ \to \pi^*(TQ_0) \to 0$ give a horizontal distribution $h_\pi = h(\pi^*(TQ_0))$ on $Q$ which is called constraint submanifold[35]. $TQ$ can be decomposed into a direct sum of horizontal and vertical space by use of projection operators[36]

$$p_h = \left(\frac{\partial}{\partial q^\mu} + B_\mu^\alpha \frac{\partial}{\partial q^\alpha}\right) \otimes dq^\mu, \quad p_v = \frac{\partial}{\partial q^\alpha} \otimes \left(dq^\alpha - B_\mu^\alpha dq^\mu\right) \tag{3.2}$$

Frobenius integrability of the constraints is determined by the curvature of the connection $h$, i.e., $R = dp_v \cdot p_h$, which is locally equivalent to

$$R_{\mu\nu}^\alpha = \frac{\partial B_\mu^\alpha}{\partial q^\nu} - \frac{\partial B_\nu^\alpha}{\partial q^\mu} \tag{3.3}$$

Denote by $\mathcal{L} = i^*L$ a regular Lagrangian on constraint submanifold $h_\pi$, where $i: h_\pi \to TQ$. The non-degenerate fundamental 2-form on $h_\pi$ can be constructed by

$$\mathcal{W} = d\left(\frac{\partial \mathcal{L}}{\partial \dot{q}^\mu} dq^\mu\right) + \frac{1}{2} R^\alpha_{\mu\nu} \, i^*\left(\frac{\partial L}{\partial q^\alpha}\right) dq^\mu \wedge dq^\nu \tag{3.4}$$

Let $Z = \dot{q}^\mu \, \partial/\partial q^\mu + B^\alpha_\mu \dot{q}^\mu \, \partial/\partial q^\alpha + f^\mu \, \partial/\partial \dot{q}^\mu$ be a dynamical vector field on $h_\pi$ and

$E_\mathcal{L} = \frac{\partial \mathcal{L}}{\partial \dot{q}^\mu} \dot{q}^\mu - \mathcal{L}$ the energy function. Then the Chaplygin's equations[10,11]

$$Z\left(\frac{\partial \mathcal{L}}{\partial \dot{q}^\mu}\right) - \frac{\partial \mathcal{L}}{\partial q^\mu} + R^\alpha_{\nu\mu} \, \dot{q}^\nu i^*\left(\frac{\partial L}{\partial \dot{q}^\alpha}\right) = 0 \tag{3.5}$$

can be geometrically represented by

$$i_Z \mathcal{W} = -dE_\mathcal{L} \tag{3.6}$$

Now we turn to discussion of an almost Poisson structure of the constraint submanifold. The Legendre transformation is utilized to define a constraint phase space $\mathcal{M} = \mathbb{F}L(h_\pi)$ and Hamiltonian $\mathcal{M} \to R: \mathcal{H} = p_\mu \dot{q}^\mu - \mathcal{L}$ where the momentum $p_\mu = \mathbb{F}L(\dot{q}^\mu)$. Hence the fundamental 2-form $\mathbb{F}L(\mathcal{W}) = \Omega$ on $\mathcal{M}$ is

$$\Omega = dp_\mu \wedge dq^\mu + \frac{1}{2} R^\alpha_{\mu\nu} \, p_\alpha dq^\mu \wedge dq^\nu \tag{3.7}$$

where the $p_\alpha$ is restricted to $\mathcal{M}$. The components of $\Omega$ constitute a non-symplectic matrix

$$\Omega = \begin{pmatrix} R^\alpha_{\mu\nu} p_\alpha & -\delta_{\mu\nu} \\ \delta_{\mu\nu} & 0 \end{pmatrix} \tag{3.8}$$

The almost Poisson tensor $\mathbf{J} = \Omega^{-1}$ is then

$$\mathbf{J} = \begin{pmatrix} 0 & \delta_{\mu\nu} \\ -\delta_{\mu\nu} & R^\alpha_{\mu\nu} p_\alpha \end{pmatrix} = \begin{pmatrix} 0 & \delta_{\mu\nu} \\ -\delta_{\mu\nu} & 0 \end{pmatrix} + \begin{pmatrix} 0 & 0 \\ 0 & R^\alpha_{\mu\nu} p_\alpha \end{pmatrix} \tag{3.9}$$

which is a sum of canonical Poisson tensor and almost Poisson one. The almost Hamiltonian vector field $X_\mathcal{H}$ of the function $\mathcal{H}$ is defined by

$$i_{X_\mathcal{H}} \Omega = -d\mathcal{H} \tag{3.10}$$

It should Pointed that the definition here is different from section 2 by a minus to be suitable to conventional representation. Then an almost Poisson bracket of two functions $f, \mathcal{H} \in \mathcal{F}(\mathcal{M})$ on $\mathcal{M}$ can be constructed by

$$\Omega(X_f, X_\mathcal{H}) = -X_\mathcal{H}(f) = -[f, \mathcal{H}] \tag{3.11}$$

which can also be taken as a definition of the almost Hamiltonian vector field.

The components of the almost Hamiltonian vector field $X_{\mathcal{H}}$ is then

$$X_{\mathcal{H}}(q^{\mu}) = \dot{q}^{\mu} = [q^{\mu}, \mathcal{H}] = \frac{\partial \mathcal{H}}{\partial p_{\mu}} \tag{3.12a}$$

$$X_{\mathcal{H}}(p_{\mu}) = \dot{p}_{\mu} = [p_{\mu}, \mathcal{H}] = -\frac{\partial \mathcal{H}}{\partial q^{\mu}} + R^{\alpha}_{\mu\nu} p_{\alpha} \frac{\partial \mathcal{H}}{\partial p_{\nu}} \tag{3.12b}$$

From Eq. (3.7) or Eq. (3.10b) we can decompose the almost bracket into

$$[p_{\mu}, \mathcal{H}] = \{p_{\mu}, \mathcal{H}\} + [\![p_{\mu}, \mathcal{H}]\!] \tag{3.13}$$

where

$$\{f, \mathcal{H}\} = \omega^{ij} \frac{\partial f}{\partial x^{i}} \frac{\partial \mathcal{H}}{\partial x^{j}} = \frac{\partial f}{\partial q^{\mu}} \frac{\partial \mathcal{H}}{\partial p_{\mu}} - \frac{\partial f}{\partial p_{\mu}} \frac{\partial \mathcal{H}}{\partial q^{\mu}} \tag{3.14a}$$

$$[\![f, \mathcal{H}]\!] = R^{\alpha}_{\mu\nu} p_{\alpha} \frac{\partial f}{\partial p_{\mu}} \frac{\partial \mathcal{H}}{\partial p_{\nu}} \tag{3.14b}$$

Notice that the almost Poisson bracket $[\![,]\!]$ is restricted to a closed half space spanned by $\{dp_{\mu}\}$ since

$$[\![p_{\mu}, p_{\nu}]\!] = R^{\alpha}_{\mu\nu} p_{\alpha} \tag{3.15}$$

Then the Eq. (3.14b) can be reformulated by

$$[\![f, \mathcal{H}]\!] = [\![p_{\mu}, p_{\nu}]\!] \frac{\partial f}{\partial p_{\mu}} \frac{\partial \mathcal{H}}{\partial p_{\nu}} \tag{3.16}$$

It can be verified from Eq. (2.13) that $[\![,]\!]$ is Poisson bracket if and only if

$$\left(R^{\alpha}_{\mu\nu} R^{\beta}_{\sigma\lambda} + R^{\alpha}_{\mu\sigma} R^{\beta}_{\lambda\nu} + R^{\alpha}_{\mu\lambda} R^{\beta}_{\nu\sigma}\right) p_{\alpha} \frac{\partial p_{\beta}}{\partial p_{\mu}} = 0 \tag{3.17}$$

For a mechanical system with metric $g_{ij} = \frac{\partial^{2} L}{\partial \dot{q}^{i} \partial \dot{q}^{j}} = \delta_{ij}$ ( $i, j = 1, 2, \cdots, n$ ), $p_{\alpha} = B^{\alpha}_{\mu} p_{\mu}$, $\frac{\partial p_{\beta}}{\partial p_{\mu}} = B^{\beta}_{\mu}$. Due to the decomposition of Eq. (3.13), the dynamical vector field can be decomposed into

$$X_{\mathcal{H}} = \mathcal{X}_{\mathcal{H}} + \mathbb{X}_{\mathcal{H}} \tag{3.18}$$

where

$$\mathcal{X}_{\mathcal{H}} = \frac{\partial \mathcal{H}}{\partial p_{\mu}} \frac{\partial}{\partial q^{\mu}} - \frac{\partial \mathcal{H}}{\partial q^{\mu}} \frac{\partial}{\partial p_{\mu}}; \quad \mathbb{X}_{\mathcal{H}} = R^{\alpha}_{\mu\nu} p_{\alpha} \frac{\partial \mathcal{H}}{\partial p_{\nu}} \frac{\partial}{\partial p_{\mu}} \qquad (3.19)$$

In the case of Eq. (3.17), $\mathbb{X}_{\mathcal{H}}$ becomes a Hamiltonian vector field which maybe easy to integrate.

If a symplectic structure exists for the nonholonomic system, so does a Poisson structure according to the theorem 2.3. This condition can be obtained from Eq. (2.19) and Eq. (3.20) as following

$$R^{\alpha}_{\nu\sigma} \frac{\partial p_{\alpha}}{\partial q^{\mu}} + R^{\alpha}_{\sigma\mu} \frac{\partial p_{\alpha}}{\partial q^{\nu}} + R^{\alpha}_{\mu\nu} \frac{\partial p_{\alpha}}{\partial q^{\sigma}} = 0 \qquad (3.20a)$$

$$R^{\alpha}_{\mu\nu} \frac{\partial p_{\alpha}}{\partial p_{\sigma}} = 0 \qquad (3.20b)$$

$$\frac{\partial R^{\alpha}_{\nu\sigma}}{\partial q^{\mu}} + \frac{\partial R^{\alpha}_{\sigma\mu}}{\partial q^{\nu}} + \frac{\partial R^{\alpha}_{\mu\nu}}{\partial q^{\sigma}} \equiv 0 \qquad (3.20c)$$

which can be verified to be same with that from Eq. (2.14) and Eq. (3.9) using $\mathbf{K}^{\mu\nu} = R^{\alpha}_{\mu\nu} p_{\alpha}$.

This result implies that the bracket $[\![,]\!]$ or tensor $\mathbf{K}^{\mu\nu} = R^{\alpha}_{\mu\nu} p_{\alpha}$ is Poisson because of theorem 2.2 or corollary 2.2. Indeed, Eq. (3.17) is a direct result of Eq. (3.20b). Thus we have

**Theorem 3.1**: The Poisson bracket of Chaplygin's nonholonomic systems with a symplectic structure can always decomposed into a sum of a canonical Poisson bracket and a non-canonical Poisson bracket.

This result is suitable to the decomposition of the Hamiltonian vector filed. However, if the bracket $[\![f,\mathcal{H}]\!]$ satisfies the condition (3.17) does not definitely lead to a Poisson structure $[f,\mathcal{H}]$, which will be illustrated in next section.

## 4. Almost Poisson structure associated with torsion of an affine space

There exist some physical systems, e.g., elementary particles moving in Riemann-Cartan spacetime, a crystal with dislocation, motion of rigid body in body-fixed coordinate system, *etc.*, whose physical property can be characterized geometrically by the torsion of a general affine metric space in which both connection and metric are independently taken as essential geometric objects. An autoparallel of such space will deviate from its geodesic unless the affine connection is symmetric and compatible with the metric[12-15]. The free motion of such systems in affine space is described by its autoparallels not by geodesic lines[15, 37-43].

Let $M$ denote the affine space with local coordinates $\{q^{\mu}\}$ ($i = 1, 2, \cdots, n$). Denote the metric of $M$ by $g_{\mu\nu}$ and Riemann connection derived from $g_{\mu\nu}$ by Christoffel symbol

$$\overline{\Gamma}^{\mu}_{\nu\sigma} = \frac{1}{2} g^{\mu\lambda} \left( \partial_{\sigma} g_{\nu\lambda} + \partial_{\nu} g_{\sigma\lambda} - \partial_{\lambda} g_{\nu\sigma} \right) \triangleq \frac{1}{2} g^{\mu\lambda} \overline{\Gamma}_{\lambda\nu\sigma} \qquad (4.1)$$

The affine connection is denoted by $\Gamma^{\mu}_{\nu\sigma}$ whose antisymmetric part $S^{\mu}_{\nu\sigma} \triangleq \frac{1}{2}\left(\Gamma^{\mu}_{\nu\sigma} - \Gamma^{\mu}_{\sigma\nu}\right)$ is called the torsion of the space $M$. Then the autoparallels and geodesics are respectively represented by

$$\ddot{q}^{\mu} + \Gamma^{\mu}_{\nu\sigma}\dot{q}^{\nu}\dot{q}^{\sigma} = 0 \tag{4.2a}$$

$$\ddot{q}^{\mu} + \overline{\Gamma}^{\mu}_{\nu\sigma}\dot{q}^{\nu}\dot{q}^{\sigma} = 0 \tag{4.2b}$$

Since Eq. (4.2a) plays an important role in analyzing the motion of some physical systems we take it as a starting point to discuss its inverse problem of calculus of variations. Let $S_{\nu\sigma\lambda} \triangleq g_{\nu\mu}S^{\mu}_{\sigma\lambda}$. Contracting a metric $g_{\lambda\mu}$ with the equation and making use of the relation[15]

$$g_{\mu\lambda}\Gamma^{\mu}_{\nu\sigma} = \overline{\Gamma}_{\lambda\nu\sigma} - 2S_{\nu\sigma\lambda} \tag{4.3}$$

we obtain

$$g_{\lambda\mu}\ddot{q}^{\mu} + \overline{\Gamma}_{\lambda\nu\sigma}\dot{q}^{\nu}\dot{q}^{\sigma} - 2g_{\nu\rho}S^{\rho}_{\sigma\lambda}\dot{q}^{\nu}\dot{q}^{\sigma} = 0 \tag{4.4}$$

If we require the symmetry of time reparametrization and general coordinate transformations for actions to establish a relativistic quantum theory for the autoparallel motion, the Lagrangian for the equation should be $L_g = \sqrt{g_{\mu\nu}\dot{q}^{\mu}\dot{q}^{\nu}}$. The Eq. (4.4) becomes

$$\left[L_g\right]_{\mu} + 2g_{\nu\rho}S^{\rho}_{\sigma\lambda}\dot{q}^{\nu}\dot{q}^{\sigma}/\sqrt{\dot{q}^2} = 0 \tag{4.5}$$

where $\left[L_g\right]_{\mu}$ is Lagrange derivative of $L_g$. However, the Helmholtz conditions[1] restrict the torsion force vanishes, which gives a trivial solution. Of course a self-adjoint genotopic transformation can be used to obtain Euler-Lagrange equations with a new Lagrangian $L_{\phi} = e^{\phi(q)}L_g$ where the integrating factor $e^{\phi(q)}$ is related with torsion by $S^{\lambda}_{\mu\nu} = \frac{1}{2}(\delta^{\lambda}_{\mu}\partial_{\nu}\phi(q) - \delta^{\lambda}_{\nu}\partial_{\mu}\phi(q))$. But the corresponding canonical Hamiltonian vanishes identically due to the above symmetry. Therefore, we have to remove the restriction on the time reparametrization symmetry for the actions and take $L = g_{\mu\nu}\dot{q}^{\mu}\dot{q}^{\nu}$ as a Lagrangian to get a regular Hamiltonian. In this selection, the Eq. (4.4) becomes

$$\left[L\right]_{\lambda} = 2S^{\rho}_{\lambda\sigma}\frac{\partial L}{\partial \dot{q}^{\rho}}\dot{q}^{\sigma} \tag{4.6}$$

By means of standard Legendre transformation, the Hamiltonian $T^*M \to R: H = p_{\mu}\dot{q}^{\mu} - L$ is defined, where $p_{\mu} = \mathbb{F}L(\dot{q}^{\mu}) = \partial L/\partial \dot{q}^{\mu}$ is the momentum. Then the Eq. (4.5) can be transformed into the Hamiltonian formulation:

$$\dot{q}^{\mu} = \frac{\partial H}{\partial p_{\mu}}, \qquad \dot{p}_{\mu} = -\frac{\partial H}{\partial q^{\mu}} - 2S^{\sigma}_{\mu\nu}p_{\sigma}\frac{\partial H}{\partial p_{\nu}} \tag{4.7}$$

Making use of a fundamental 2-form on $T^*M$

$$\Omega = dp_\mu \wedge dq^\mu - S^\sigma_{\mu\nu} p_\sigma dq^\mu \wedge dq^\nu \tag{4.8}$$

the Eq. (4.7) can be formulated geometrically by

$$i_{X_H}\Omega = -dH \tag{4.9}$$

which can also be taken as a definition of an almost Hamiltonian vector field $X_H$ on $T^*M$.

Define an almost Poisson bracket operation $[,]$ on the set of functions $\mathcal{F}(T^*M)$ on $T^*M$ by

$$i_{X_f}i_{X_H}\Omega = X_H(f) = [f, H] \tag{4.10}$$

Due to the torsion, the bracket is an almost Poisson one and accordingly the pair $\{T^*M,[,]\}$ is an almost Poisson manifold. This bracket can be decomposed into a sum of canonical Poisson bracket and an almost Poisson one

$$[f, H] = \{f, H\} + [\![f, H]\!]$$

where

$$\{f, H\} = \frac{\partial f}{\partial q^\mu}\frac{\partial H}{\partial p_\mu} - \frac{\partial f}{\partial p_\mu}\frac{\partial H}{\partial q^\mu}, \quad [\![f, H]\!] = 2S^\sigma_{\mu\nu} p_\sigma \frac{\partial f}{\partial p_\mu}\frac{\partial H}{\partial p_\nu} \tag{4.11}$$

In terms of $\{p_\mu\}$ the almost Poisson bracket reduces to

$$[\![p_\mu, p_\nu]\!] = 2S^\sigma_{\mu\nu} p_\sigma \tag{4.12}$$

which implies that the half space spanned by $\{dp_\mu\}$ is closed. If torsion tensor $S^\sigma_{\mu\nu}$ satisfies the following condition

$$S^\sigma_{\mu\nu}S^\mu_{\rho\lambda} + S^\sigma_{\mu\rho}S^\mu_{\lambda\nu} + S^\sigma_{\mu\lambda}S^\mu_{\nu\rho} = 0 \tag{4.13}$$

it will play a role on the half space $T^*M/M$ as structure constants of a Lie group. Then the almost Poisson bracket $[\![f, H]\!]$ becomes a "Lie-Poisson bracket". Otherwise it is an almost Lie-Poisson bracket. With the help of Eq. (4.12) the almost Poisson bracket $[\![f, H]\!]$ can also be represented by

$$[\![f, H]\!] = [\![p_\mu, p_\nu]\!]\frac{\partial f}{\partial p_\mu}\frac{\partial H}{\partial p_\nu} \tag{4.14}$$

By means of the brackets the Eq.(4.7) becomes

$$\dot{q}^\mu = X_H(q^\mu) = \{q^\mu, H\}, \quad \dot{p}_\mu = X_H(p_\mu) = \{p_\mu, H\} + [\![p_\mu, H]\!] \tag{4.15}$$

Accordingly the almost Hamiltonian vector field $X_H$ can be decomposed into a sum of a canonical Hamiltonian vector field and an almost Hamiltonian vector field

$$X_H = \mathcal{X}_H + \mathbb{X}_H \tag{4.16}$$

where

$$\mathcal{X}_H = \frac{\partial H}{\partial p_\mu}\frac{\partial}{\partial q^\mu} - \frac{\partial H}{\partial q^\mu}\frac{\partial}{\partial p_\mu}; \quad \mathbb{X}_{\mathcal{H}} = -2S^\alpha_{\mu\nu} p_\alpha \frac{\partial H}{\partial p_\nu}\frac{\partial}{\partial p_\mu} \tag{4.17}$$

## 5. Applications to a homogeneous ball rolling without slipping

Let us discuss the almost Poisson structure on the constraint submanifold of a homogeneous ball with mass $m$ and radius $a$ freely rolling on a rough horizontal plane. The configuration $Q$ of the rolling ball is recognized by coordinates $(x, y)$ of the centre of mass of the ball and the three Euler angles $(\psi, \theta, \varphi)$. The non-holonomic constraints are

$$\dot{x} = -a(\dot{\varphi}\sin\theta\cos\psi - \dot{\theta}\sin\psi), \quad \dot{y} = -a(\dot{\varphi}\sin\theta\sin\psi + \dot{\theta}\cos\psi)$$

With the notation of section 3, we obtain

$$B^4_1 = 0, \; B^4_2 = \sin\psi, \; B^4_3 = -a\sin\theta\cos\psi; \; B^5_1 = 0, \; B^5_2 = -\cos\psi, \; B^5_3 = -a\sin\theta\sin\psi$$

The Lagrangian of the system is $L = \frac{1}{2}m(\dot{x}^2 + \dot{y}^2) + \frac{1}{2}\cdot\frac{2}{5}ma^2(\dot{\psi}^2 + \dot{\theta}^2 + \dot{\varphi}^2 + 2\dot{\psi}\dot{\varphi}\cos\theta)$.

Taking notational identifications: $q^\mu = (\psi, \theta, \varphi)$ $(\mu = 1, 2, 3)$, $q^\alpha = (x, y)$ $(\alpha = 4, 5)$, we have

$$R^4_{12} = a\cos\psi, \; R^4_{13} = a\sin\theta\sin\psi, \; R^4_{23} = -a\cos\theta\cos\psi$$

$$R^5_{12} = -a\sin\psi, \; R^5_{13} = -a\sin\theta\cos\psi, \; R^5_{23} = -a\cos\theta\sin\psi$$

Then the Lagrangian pull-back to the constraint submanifold $h_\pi$ is

$$\mathcal{L} = \frac{1}{2}ma^2\left[\frac{7}{5}\dot{\theta}^2 + \dot{\varphi}^2\sin^2\theta + \frac{2}{5}(\dot{\psi}^2 + \dot{\varphi}^2 + 2\dot{\psi}\dot{\varphi}\cos\theta)\right]$$

Making use of Legendre transformation, the corresponding Hamiltonian on the submanifold $\mathcal{M} = \mathbb{F}L(h_\pi)$ is given by

$$\mathcal{H} = \frac{5}{2ma^2}\left[\frac{5p_\psi^2}{4} + \frac{p_\theta^2}{7} - \frac{5\cos 2\theta (p_\psi - p_\varphi)^2}{49\sin^4\theta}\right]$$

After complicated but direct computation the almost Poisson tensor can be obtained with a decomposition

$$\mathbf{J} = \begin{pmatrix} 0 & 0 & 0 & 1 & 0 & 0 \\ 0 & 0 & 0 & 0 & 1 & 0 \\ 0 & 0 & 0 & 0 & 0 & 1 \\ -1 & 0 & 0 & 0 & 0 & 0 \\ 0 & -1 & 0 & 0 & 0 & 0 \\ 0 & 0 & -1 & 0 & 0 & 0 \end{pmatrix} + \begin{pmatrix} 0 & 0 & 0 & 0 & 0 & 0 \\ 0 & 0 & 0 & 0 & 0 & 0 \\ 0 & 0 & 0 & 0 & 0 & 0 \\ 0 & 0 & 0 & 0 & \mathbf{K}_{12} & \mathbf{K}_{13} \\ 0 & 0 & 0 & -\mathbf{K}_{12} & 0 & \mathbf{K}_{23} \\ 0 & 0 & 0 & -\mathbf{K}_{13} & -\mathbf{K}_{23} & 0 \end{pmatrix}$$

where

$$\mathbf{K}_{12} = \frac{5}{7}\left(p_\theta \sin 2\psi - \frac{\cos 2\psi}{\sin\theta}(p_\psi - p_\varphi)\right),\ \mathbf{K}_{13} = \frac{5}{7} p_\theta \sin\theta,\ \mathbf{K}_{23} = \frac{5}{7}\cot\theta (p_\psi - p_\varphi)$$

Based on this tensor the almost Poisson bracket and almost Hamiltonian vector field can be represented by the former standard procedure given in section 3. For example, the equations of motion can be formulated by

$$\dot{\psi} = \{\psi, \mathcal{H}\},\ \dot{\theta} = \{\theta, \mathcal{H}\},\ \dot{\varphi} = \{\varphi, \mathcal{H}\}$$

$$\dot{p}_\psi = \{p_\psi, \mathcal{H}\} + [\![p_\psi, \mathcal{H}]\!],\ \dot{p}_\theta = \{p_\theta, \mathcal{H}\} + [\![p_\theta, \mathcal{H}]\!],\ \dot{p}_\varphi = \{p_\varphi, \mathcal{H}\} + [\![p_\varphi, \mathcal{H}]\!]$$

where $[\![p_\mu, \mathcal{H}]\!] = \mathbf{K}_{\mu\nu}\dfrac{\partial \mathcal{H}}{\partial p_\nu}$ $(\mu, \nu = 1, 2, 3)$ is an almost Poisson bracket on the momentum space.

It should be pointed that the configuration of a Chaplygin's nonholonomic system admit a affine space structure with torsion defined by[15]

$$S^\rho_{\mu\nu} = \frac{1}{2} g^{\rho\sigma} g_{ij} \varepsilon^i_\sigma \left(\frac{\partial \varepsilon^j_\nu}{\partial q^\mu} - \frac{\partial \varepsilon^j_\mu}{\partial q^\nu}\right)\ (i, j = 1, 2, 3, 4, 5)$$

where $\varepsilon^\alpha_\rho = B^\alpha_\rho, \varepsilon^\sigma_\rho = \delta^\sigma_\rho$, $g_{ij} = \partial^2 L / \partial \dot{q}^i \partial \dot{q}^j$ is an element of Hessian matrix and $g^{\rho\sigma}$ is the inverse of metric $g_{\rho\sigma} = g_{ij}\varepsilon^i_\rho \varepsilon^j_\sigma$. The autoparallels of the space are integral curves of Chaplygin's equations. The torsion tensor $S^\rho_{\mu\nu}$ is proportional to Poisson tensor $\mathbf{K}_{\mu\nu}$. Thus this example can also be utilized to construct an almost Poisson structure associated with torsion of affine space. It is interesting that if we choose an alternative Lagrangian or Hamiltonian, a symplectic or Poisson structure can be indirectly constructed[35].

# 6. Concluding Remarks

A theory of almost Poisson structure for non-self-adjoint dynamical systems is formulated in the direction of decomposition of an almost Poisson bracket into a sum of Poisson one and almost Poisson one. In many mechanical systems the latter almost Poisson bracket become almost Lie-Poisson one on a half space cotangent to the configuration manifold. In order to integrate dynamical vector field easily we further hope the almost Lie-Poisson bracket to be a real Lie-Poisson one for some mechanicall systems, which maybe useful to Poisson reduction of mechanical systems with symmetries. A good candidate of such systems is given in Ref. 44 and other candidates may appear in nonholonomic systems on Lie groups or Lagrangian mechanics or Hamiltonian mechanics on Lie algebroids, which will be investigated in the succeeding article. After the decomposition and classification of the almost Poisson structure it is still an important but difficult task to find its integrability.

**Acknowledgement** This research was partially supported by the National Nature Science Foundational of China (Grant Nos. 10872084,10472040), the Outstanding Young Talents Training Fund of Liaoning Province, China (Grant No.3040005) and the Research Program of Higher Education of Liaoning Province, China (Gant No. 2008S098).